%% file: cocentral_kac.tex
\documentclass[12pt]{amsart}
\usepackage{amsmath,amsthm,amsfonts,amssymb}
\usepackage{amsmath, amsthm, amssymb, amscd, amsfonts}
\usepackage[all]{xy}
\usepackage{amssymb, amsthm, amsmath}
\usepackage[english]{babel}
\usepackage{amsfonts}
\usepackage{color}

\newcommand{\ra}{\rightarrow}

\newcommand{\ot}{\otimes}
\newcommand{\mtc}{\mathcal}
\newcommand{\lam}{\lambda}
\newcommand{\lb}{\label}
\newcommand{\Lam}{\Lambda}

\newcommand{\al}{\alpha}
\newcommand{\ovl}{\overline}
\newcommand{\eps}{\epsilon}
\newcommand{\bn}{\begin}

\newcommand{\krn}{\mtr{ker}}
\newcommand{\Irr}{\mtr{Irr}}
\newcommand{\cS}{\mathcal{S}}\newcommand{\Rep}{\mtr{Rep}}
\newcommand{\cO}{\mtc{O}}\newcommand{\X}{\mtc{X}}
\newcommand{\B}{\mtc{B}}

\newcommand{\D}{\mtc{D}}

\newcommand{\onh}{On the other hand}

\newcommand{\rh}{\rightharpoonup}
\newcommand{\lh}{\leftharpoonup}
\numberwithin{equation}{section}

\newtheorem{defn}[equation]{Definition}
\newtheorem{cor}[equation]{Corollary}
\newtheorem{rem}[equation]{Remark}

\newcommand{\bp}{\bn{prop}}
\newcommand{\ep}{\end{prop}}

\newcommand{\dw}{\downarrow}
\newcommand{\uw}{\uparrow}

\newcommand{\ch}{\chi}
\newcommand{\mtr}{\mathrm}

\newcommand{\ncm}{\newcommand}\newcommand{\gm}{\gamma}
\numberwithin{equation}{section}
\newcommand{\el}{\end{lemma}}\newcommand{\bl}{\bn{lemma}}
\newcommand{\et}{\end{thm}}\newcommand{\bt}{\bn{thm}}
\newcommand{\beqarn}{\begin{eqnarray*}}
\newcommand{\eeqarn}{\end{eqnarray*}}
\newcommand{\beqn}{\bn{equation*}}
\newcommand{\eeqn}{\end{equation*}}
\newcommand{\bpf}{\bn{proof}}
\newcommand{\epf}{\end{proof}}\newcommand{\di}{\delta}
\ncm{\cX}{\mtc{X}}
\ncm{\wt}{\widetilde}
\ncm{\sg}{\sigma}
\ncm{\np}{\newpage}
\ncm{\ebl}{\end{thebibliography}}
\ncm{\bbl}{\begin{thebibliography}}
\ncm{\chd}{_{ _{\ch}}}
\ncm{\ald}{_{ _{\al}}}
\ncm{\cP}{\mathcal{P}}
\ncm{\ei}{e_i}
\ncm{\eij}{e_{i,\;j}}
\ncm{\bdef}{\begin{defn}}
\ncm{\edf}{\end{defn}}
\ncm{\bc}{\begin{cor}}
\ncm{\ec}{\end{cor}}
\ncm{\er}{\end{rem}}
\ncm{\br}{\begin{rem}}
\ncm{\bd}{\begin{document}}
\ncm{\ed}{\end{document}}
\ncm{\beq}{\begin{equation}}
\ncm{\dzs}{\#_{\sg}}

\ncm{\bne}{\bn{enumerate}}
\ncm{\ene}{\end{enumerate}}
\ncm{\ovr}{\overline}
\ncm{\blue}{\textcolor[rgb]{.00, .00, 1.00}}
\ncm{\red}{\textcolor[rgb]{1.00, .00, .00}}
\ncm{\green}{\textcolor[rgb]{.00, 1.00, .00}}
\ncm{\md}{\medbreak}
\input{macros}
\begin{document}
\title[Semisimple Hopf algebras]
{On coideal subalgebras of cocentral Kac algebras and a generalization of Wall's conjecture}
\begin{abstract}
It shown that any coideal subalgebra of a finite dimensional Hopf algebra is a cyclic module over the dual Hopf algebra.  Using this we describe all coideal subalgebras of a cocentral abelian extension of Hopf algebras extending some results from \cite{Grln}. %In particular a conjecture proposed by R. Guralnick for $"^*"$- Kac algebras also holds for Hopf Kac alegbras.
\end{abstract}
%They are all group theoretical. 
\author{Sebastian  Burciu}
\address{Inst.\ of Math.\ ``Simion Stoilow" of the Romanian Academy
P.O. Box 1-764, RO-014700, Bucharest, Romania}\address{and}
\address{University of Bucharest, Faculty of Mathematics and Computer Science, Algebra and Number Theory Research Center, 14 Academiei St., Bucharest, Romania }\email{sebastian.burciu@imar.ro} 
\thanks{This work was supported by a grant of the Romanian National Authority for Scientific Research, CNCS-UEFISCDI, grant no.
88/05.10.2011.}%{This research was partially supported by CNCSIS PN II RU-PD grant "Representations of semisimple Hopf algebras and fusion categories" PD-168 no: 14/28.07.2010}
%\begin{document}
\subjclass[2000]{Primary 16W30, 18D10}
%\keywords{Double crossed products; Drinfeld doubles; Fusion categories; Clifford theory}
\keywords{cocentral abelian extensions; coideal subalgebras}
 \maketitle
%\tableofcontents
\section{Introduction}
%importance of Kac algebras
%importance of fusion categories ; a systematic study DGNO, ENO
Kac algebras were among the first given examples of noncommutative noncocommutative Hopf algebras \cite{kac-p}. They are also called abelian extensions since they satisfy the following long exact sequence of Hopf algebras:
\beq\lb{defi}
k\ra k^G \ra A\ra kF\ra k
\eeq
where $G$ and $F$ are finite groups and $k^G$ is the dual Hopf algebra of the group algebra $kG$.
Irreducible representations of abelian extensions are completely characterized in \cite{KMM}.

%Twisted Drinfeld doubles of finite groups are examples of Kac algebras. Here the tensor product of two such representations is well known, \cite{cibils} and \cite{scoh}. In order to described all the braided group theoretical fusion categories, recently the authors of \cite{NNW} classify all the fusion subcategories of the category of finite dimensional representations of a twisted Drinfeld double $D^{\omega}(G)$ of a finite group $G$. Their classification is parameterized by a pair of commuting subgroups $H$ and $K$ of $G$ and a linear bicharacter $B:K\times H\ra k^*$ defined on their product.
\medbreak
\ncm{\rep}{\Rep}
Our main result describes all left (right) coideal subalgebras of cocentral Kac algebras:
\bn{thm}\label{main}
Let $A=k^G\;^{\tau}\#_{\sg} kF$ be an abelian cocentral extension of $kF$ by $k^G$. Then left coideal subalgebras of $A$ are parameterized by the following data:
\bne
\item Two subgroups $M\leq G$ and $H\leq F$ with $H\rhd M=M$ 
\item A twisted bicharacter $\lam: M\times H\ra U(1)$ satisfying the following properties
\begin{equation}\label{comp1}
    \lam(ab,\;h)=\lam(a,\;h) \lam(b,\;h)\tau_h(a,\;b)
\end{equation}

\begin{equation}\label{comp2}
   \lam(a,\;hl)= \lam(a,\;h)\lam(h _1 \rhd a,\;l)\sigma_a(h,\;l)
\end{equation}
\ene
The coideal subalgebra corresponding tothe triple $(M,H,\lam)$ is given by $$C(M,H,\;\lam):=\bigoplus_{h \in H}C_{\lam}(h)\dzs h,$$
where the space $C_{\lam}(h)$ is defined by 
$$C_{\lam}(h):=\{f \in k^G\;|R_{m, \tau_h}(f)=\lam(m,h)f \;\text{for all}\;m \in M\}.$$
\end{thm}
Here $\rhd$ is the induced action of $F$ on $G$. For the definition of the operator $R_{m, \tau_h}$ see Equation \ref{r}. 
\md
In \cite{Grln} the authors described all left (right) coideal subalgebras of Kac algebras of Izumi-Kosaki type. Using this they proved one of the Wall's conjecture 
for these algebras in the case that both finite groups $F$ and $G$ are solvable. These Kac algebras are introduced in \cite{14} and are Hopf algebras with an additional $\mathbb{C}^*$-structure. They are also studied in more details
in \cite{IK} by considering compositions of group type subfactors. 
\md
In  this paper we show that one has the same description of the coideal subalgebras as in \cite{Grln} even without the assumption on the presence of the additional $\mathbb{C}^*$-structure. This is compensated by the new characterization of coideal subalgebras given in Theorem \ref{dec}.  Also a Hopf algebraic version of the Conjecture 1.1 formulated in \cite{Grln} follows for cocentral Kac algebras of solvable groups, see Subsection \ref{wall}. 

\md
%In this paper we describe all the fusion subcategories of the category of finite dimensional representations of a cocentral finite dimensional Kac algebra $A$. It is known (see \cite{N1}) that $\Rep(A)$ is a group theoretical category in the sense defined in \cite{ENO}. Since the category $\rep(A)$ is not braided in general one cannot use Mueger's centralizer of a fusion subcategory  defined only for braided fusion categories.% in \cite{Mu}. 
%Thus our approach is different from the one used in \cite{NNW} and it is based on a systematic use of the tensor product formula for two simple $A$-modules given in \cite{scoh}.

%On the other hand we also give a conceptual explanation of the result. It is recently shown in \cite{Grln} that all coideal subalgebras of $A$ are of parameterized by a similar data. Two subgroups $M$ of $G$ and $H$ of $F$ and a twisted bicharacter $\lam:M\times H\ra k^*$.
%\blue{Natale, banica adv math; coideal subalgebras important}
Note that recently, in \cite{bbn} the authors proved that a cosemisimple Hopf algebra $A$ is a quantum permutation algebra if and only if it is generated as an algebra by the matrix coefficients of all its left (right) coideal subalgebras of $A$. %Knowing all left coideal subalgebras of a given Hopf algebras We also noticed an important restl
\md

It is well known that any finite dimensional Hopf algebra $A$ is a cyclic right $A^*$-module generated by the left integral of $\Lam \in A$. In Theorem \ref{dec} we prove an analogue of this result for left (right) coideal subalgebras of $A$. It is shown that any left coideal subalgebra $S$ of $A$ is a cyclic $A^*$-module generated by an invariant element $y_S \in S$ introduced by the author in \cite{gmj}.
%\blue{An analogue of integral, developed in \cite{glasgow} and further studied here in Theorem 2.1}
%{\bf In Theorem \ref{dec} we prove a new result on the structure of left coidela subalgebras that will enable us to extend the results from \cite{grln} by dropping the additional $*$ structure on $A$.}
\blue{
%We also show that abelian cocentral extensions all have the property $(N)$ defined in \cite{kernn}. Recall that a semisimple Hopf algebra is said to have property $(N)$ if the kernel (as defined in \cite{Bker}) of any irreducible representation of $A$ is a normal Hopf subalgebra of $A$. As shown in \cite{kernn} that is equivalent that all Hopf projections map starting from $A$ have normal Hopf kernels. It was previously shown that all Drinfeld doubles $D(G)$ of finite groups have property $(N)$ (see \cite{dakernn}).
}
%Characterization of Hopf algebras with $(N)$ given in \cite{kernn}
%more introduction from sonia, nikshych my paper
% Wall's conjecture ostrik
\md
Shortly, the organization of the paper is as follows. In the second Section we prove the result concerning the structure of coideal subalgebras as cyclic modules over the dual Hopf algebra. Section 3 contains the proof of Theorem \ref{main}. In the last section using the results from the previous section we also describe all Hopf subalgebras of cocentral abelian extensions of Hopf algebras.
\md
We work over an algebraically closed field of arbitrary characteristic and all the other Hopf algebra notations are those used in \cite{Montg}.
\section{ Preliminaries}
%\subsection{On left coideal subalgebras}
Recall that a left coideal subalgebra $S$ of $H$ is a subalgebra $S$ of $H$ with $\Delta(S)\subset H \ot S$. A coideal subalgebra $S$ of $H$ is called normal if $S$ is stable under the adjoint action of $H$ on itself, i.e $h_1sS(h_2)\in S$ for all $s \in S$ and $h \in H$.
 
%\blue{ For an arbitrary subspace $S$ of a Hopf algebra $H$ we denote by $\omega(S)$ the vector space generated by $s-\eps(s)1$ for all $s \in S$.}
%\subsection{The correspondence between Hopf ideals and normal left coideal subalgebras}\lb{corresp}%where need it?
% ON WHAT FIELD THIS WORK? FOR WHAT TYPE OF HOPF ALGEBRAS? FINITE DIMENSIONAL POINTED? ETC
%The correspondence between Hopf ideals and normal coideal subalgebras from \cite{Tkq} can be summarized as follows:

\subsection{Invariant elements of coideal subalgebras}\label{integral}
Let $S$ be any right coideal subalgebra of a finite dimensional Hopf algebra $A$. Then $A$ is free over $S$ \cite{Sk} both as left and right $S$-module. Let $A=S\oplus(\oplus_{i=2}^rSx_i)$ be a decomposition of $A$ as a free left $S$-modules of rank one. Then the idempotent integral $\Lam$ of $A$ admits a decomposition: $$\Lam=y_S+\sum_{i=2}^ry_i$$ with $y_S \in S$ and $y_i \in Sx_i$ for $i \geq 2$. Then equation $s\Lam=\eps(s)\Lam$ implies that $sy_S=\epsilon(s)y_S$ for all $s \in S$. 
\md
%If $A$ is semisimple note that $\eps(y_1)\neq 0$ by Theorem of \cite{gmj}. 
\md
Recall the left and right action of $A^*$ on $A$ given by $f\rh a=f(a_2)a_1$ and $a\lh f=f(a_1)a_2$. Next Theorem shows that any coideal subalgebra $S$ of $A$ is a cyclic right $A^*$-submodule of $A$ generated $y_S\in S$.

\bt\label{dec}
Let $S$ be a right coideal subalgebra of $A$ and $y_S\in S$ be defined as above. Then $y_S \lh A^*=S$.
\et
\bpf
From Theorem 6.1 of \cite{Sk} one has that $S$ is a simple object of the category $\;_{S}\mtc{M}^A$ of relative modules. Clearly $y_S \lh A^*\subset S$. We will show that $y_S\lh A^*\in \;_{S}\mtc{M}^A$ and then the proof will be complete.

Since $ry_S=\eps(r)y_S$ for all $r\in S$ by applying $\Delta$ it follows that $$\sum r_1(y_S)_1\ot r_2(y_S)_2=\eps(r)(y_S)_1\ot (y_S)_2.$$

Thus
$$\sum (y_S)_1\ot r(y_S)_2=S(r_1)(r_2)_1(y_S)_1\ot(r_2)_2 (y_S)_2=S(r)(y_S)_1\ot (y_S)_2,$$
for all $r\in S$. Thus $r(y_S\lh f)=f((y_S)_1)r(y_S)_2=f(S(r)(y_S)_1)(y_S)_2\in y_S\lh A^*$. This shows that $y_S \lh A^*$ is an $S$-module and therefore an object in the category $\;_S\mtc{M}^A$.
\epf

\br \lb{dualcom} Suppose that $A$ is a finite dimensional Hopf algebra and let $\{e_i\}_{i=1}^s$ be a basis of $A$ and $\{e_i^*\}_{i=1}^s \in A^*$ be the dual basis. Then for any $f \in A^*$ one has: 
$$\Delta_{A^*}(f)=\sum_{i=1}^se_i^*\otimes f \lh e_i=\sum_{i=1}^se_i \rh f\otimes e_i^*$$
\er

Indeed $f(xy)=\sum_{i=1}^se_i^*(x)f (e_iy)=\sum_{i=1}^se_i^*(x)(f \lh e_i)(y)$ for all $x,y \in B$. This shows the first equality and the second is proven similarly.

\subsubsection{On the operators $L$ and $R$ of $k^G$} Let $G$ be a finite group and $k^G$ be the dual group algebra. Consider the following operators on $k^G$ given by
$R_m(f)=f \leftharpoonup m$ and $L_m(f)= m\rightharpoonup f$ for all $m \in G$.

If $M$ is a subgroup of $G$ let $A=k^{(G/M)_l}$ be the space of all linear functions on $G$ which are constant on the right cosets of $M$ in $G$. Thus 
\beq
k^{(G/M)_l}=\{f \in k^G\;|\; f(gm)=f(m)\;\text{for all} \; g\in G \; \text{and} \;m \in M\}
\eeq

Thus $k^{(G/M)_l}$ is the subspace of all functionals $f \in k^G$ such that $L_m(f)=f$ for all $m \in M$.

\bn{lem}\label{existence}
Suppose that $A \subset k^G$ is a subalgebra of $k^G$ such that $L_g(A)=A$ for all $g \in G$. Then there is a subgroup $M$ of $G$ such that $A=k^{(G/M)_l}$.
\end{lem}

\bpf Since $L_g(A)=A$ for all $g \in G$ it follows from Formula \ref{dualcom} that $\Delta_{k^G}(A)\subset k^G \ot A$.
Thus $A$ is a normal left coideal subalgebra of $k^G$ since $k^G$ is commutative. Then it follows that the quotient $(k^G//A)^*$ is a Hopf subalgebra of $kG$. Therefore there is a subgroup $M$ of $G$ such that $ (k^G//A)^*\cong kM$. This implies that $A=k^{(G/M)_l}$.
\epf
 %%%%%%%%%%%%%%%%%%%%%%%%%%%%%%%%%%%%%%%%%%%

\section{Structure of coideal subalgebras of cocentral abelian extensions of Hopf algebras }
\subsection{The Hopf algebra $A$.}Let $A=k^G\;^{\tau}\dzs  kF$ be an arbitrary cocentral abelian extension of $k^G$ via $kF$.
Recall \cite{Nr} that this means that $A$ fits into the following exact sequence of Hopf algebras: $$k \ra k^G \ra A \ra kF \ra k.$$
Moreover the group $F$ acts by automorphisms via $\rhd : F\times G \ra G$ on the group
$G$. This induces an action of $F$ on the dual Hopf algebra $k^G$ via $f  \rh p_a=p_{f \rhd a}$.

Then the Hopf algebra $A$ has the following multiplication structure
\begin{equation}\lb{mms}
    (p_a\dzs h)(p_b\dzs l)=\delta_{a, h \rhd b}\sigma_a(h, l)p_a\dzs hl
\end{equation}
and the comultiplication structure given by
\begin{equation}\lb{comms}
  \Delta(p_a \dzs h)=\sum_{b \in G}\tau_{a}(b,\;b^{-1}a)(p_b\dzs h )\ot(p_{b^{-1}a}\dzs h).
\end{equation}
\md
 Here $\sigma: F\times F \ra k^G$ is a normalized twisted two cocycle of $G$ with respect to the action of $F$, i.e $\sg$ satisfies:
\begin{equation}\label{sg1}
    \sigma_a(h,\;l)\sigma_a(hl,\;t)=\sigma_a(h,\;lt)\sigma_{h \rhd a}(l,\;t)
\end{equation}
where by definition $\sg(h,l):=\sum_{a \in G}\sg_a(h,l)p_a$.

The dual cocycle $\tau: F \ra k^G\ot k^G$ it is denoted by $$\tau(f)=\sum_{a,b\in G}\tau_f(a,b)p_a\ot p_b$$ and satisfies the following two cocycle property:

\beq\lb{tau}
\tau_f(ab,c)\tau_f(a,b)=\tau_f(a,bc)\tau_f(b,c)
\eeq
for all $f \in F$ and $a,b,c \in G$.

Moreover one of the compatibility conditions between $\sg$ and $\tau$ is called Pentagon Equation \cite{AD} and it can be written as 
\begin{equation}\label{sg2-pent}
  \frac{\sigma_{ab}(h,\;l)}{\sigma_a(h,\;l)\sigma_b(h,\;l)}=\frac{\tau_{h}(a,\;b)\tau_l(h\rhd a,\;h \rhd b)}{\tau_{hl}(a,b)}.
\end{equation}
We can also assume that both $\sg$ and $\tau$ are normalized coccycles, that is:
\beq
\sg_g(1,f)=\sg_g(f,1)=\tau_f(1,g)=\tau_f(g,1)=1
\eeq
for all $g \in G$ and $f \in F$.
The antipode of $A$ is given by
\begin{eqnarray}\lb{ant}
    S(p_a\dzs \ovr{h}) & = & \sigma^{-1}_{h^{-1}\rhd a^{-1}}(h^{-1}, h)\tau^{-1} _{a^{-1},\;a}(h)p_{h^{-1}\rhd a^{-1}}\dzs\overline{h^{-1}}
%\\ & = & \sigma^{-1}_{h^{-1}\rhd a^{-1}}(h^{-1}, h)\tau^{-1} _{a^{-1},\;a}(h)\overline{h^{-1}}p_{a^{-1}}
\end{eqnarray}
\subsection{Left coideal subalgebras of $A$}\label{Guralnick}
In this subsection we give a complete description of all left coideal subalgebras of $A$. 
\ncm{\bq}{\beq}\ncm{\eq}{\eeq}
\subsection{Operators $L_{m ,\tau_h}$ and $R_{m ,\tau_h}$}
Define the linear operators on $k^G$ by
\bq\lb{l}
(L_{a,\tau_h}(f))=(a \rh f)\tau_h(-,\;a )
\eq
and
\bq\lb{r}
(R_{a,\tau_h}(f))=(f \lh a)\tau_h(a,-)
\eq	 
for all $a\in G$, $h \in  F$. Note that by their definition $L_{a,\tau_h}$ and $R_{a,\tau_h}$ satisfy
$$L_{a, \tau_h}(f)(b) := f(ba)\tau_h(b,a),\;\;\text{and}\;\; R_{a,\tau_h}(f)(b) := f(ba)\tau_h(a,b),$$ for all $a,b \in G$ and $h \in  F$.

The following lemma summarize the properties of these operators which follow
from definitions:
\bn{lemma}With the above notations one has that
\begin{equation}\lb{leff}
L_{a,\tau_h}L_{b,\tau_h} = L_{ab,\tau_h}\tau_h(a, b)\;\;\;\text{and}\;\:\;R_{a,\tau_h}R_{b,\tau_h} = R_{ab,\tau_h}\tau_h(b, a)
\end{equation}
Also,
\begin{equation}\lb{comlr}
   L_{a,\tau_h}R_{b,\tau_h} = R_{b,\tau_h}L_{a,\tau_h}.
\end{equation}
for all $a,b \in G$and all $h \in F$.
\end{lemma}
\bpf
The proof is by a straightforward computation.
\epf

\bl\lb{commf} Let $A=k^G\;^{\tau}\dzs  kF$ be a cocentral abelian extension of $k^G$ via $kF$.
Then the comultiplication in $A$ is given by
\beq
\Delta(u\dzs h)=\sum_{a\in G}(L_{a,\tau_h}(u)\dzs h)\ot (p_a\dzs h)=\sum_{a\in G}(p_a\dzs h)\ot (R_{a,\tau_h}(u)\dzs h)
\eeq
for any $u \in k^G$ and any $h \in F$.
\el

\bpf
We prove the first formula. Using Formula \ref{comms} one has
\beqarn
\Delta(u\dzs h) & = & \sum_{a, b \in G}\tau_h(a, b)(u_1p_a\dzs h)\ot (u_2p_b\dzs h)\\ &= & \sum_{b \in G}(\sum_{a \in G}\tau_h(a, b)u(ab)p_a\dzs h)\ot (p_b\dzs h)\\ &= & \sum_{b \in G}(L_{b,\tau_h}(u)\dzs h)\ot (p_b\dzs h).
\eeqarn
The second formula has a similar proof.
\epf

\subsection{Definition of the left coideal subalgebra $C(M,\;H,\;,\lam)$}
Let $M \leq G$ be a subgroup of $G$ and  $H \leq F$ be a subgroup of $F$ such that $M$ is stable under the action of $H$ on $G$, i.e $H \rhd M=M$. Let also  $\lam:M\times H\ra k^*$ be a twisted bicharacter on $M \times H$, i.e a function satisfying the following properties:
\beq\lb{c1}
\lam(mn,\; h)=\lam(m,\;h)\lam(n,\;h)\tau_h^{-1}(m,n)
\eeq
\beq\lb{c2}
\lam(m,\;hl)=\lam(m,\;h)\lam(h\rhd m ,\;l)\sg_{m}(h, l)
\eeq
for all $m, n\in M$ and $h, l \in H$.

Define the following subspace of $A$:
\beq
C(M,\;H,\;\lam)=\bigoplus_{h \in H}C_{\lam}(h)\dzs h
\eeq
where $$C_{\lam}(h)=\{f \in k^G\;|L_{m, \tau_h}(f)=\lam(m,h)f \;\text{for all}\;m \in M\}.$$
%%%%%%%%%%%%%%%%%%%%%%%%%%%%%%%%%%%%
\subsubsection{Description of $C_{\lam}(h)$}
\bl\lb{desciptch}
Let $f \in k^G$. One has that $f \in C_{\lam}(h)$ if and only if 
\bq\lb{gm}
f(gm)=\frac{\lam(m,\;h)}{\tau_h(g,m)}f(g)
\eq 
for all $m \in M$ and $g \in G$. In particular $\dim_k C_{\lam}(h)=\frac{|G|}{|M|}$.
\el

\bpf
Note that 
\beqarn L_{m,\tau_h}(f) & = & (m \rightharpoonup f)\; \tau_h(-,\;m)\\ & = &( m \rightharpoonup (\sum_{g\in G}f(g)p_g) )\tau_h(-,\;m)\\ & = & (\sum_{g\in G}f(g)p_{gm^{-1}})\tau_h(-,\;m)\\ & = & \sum_{g \in G}f(g)\tau_h(gm^{-1},\;m)p_{gm^{-1}}\\ & = & \sum_{g \in G}f(gm)\tau_h(g,\;m)p_{g}.\eeqarn Thus $L_{m,\tau_h}(f)=\lam(m,\;h)f$ if and only if 
\beq
f(gm)=\frac{\lam(m,\;h)}{\tau_h(g,m)}f(g)
\eeq

for all $g \in G$.

%right cosets
Let $b_i$ be a set of right coset representatives of $M$ in $G$. Thus one has $G=\sqcup_{i=1}^sb_iM$. For any $g \in G$ define the function 
\beq\lb{basis}
f_{[g]}=\sum_{m \in M}\frac{\lam(m,\;h)}{\tau_h(g,m)}p_{gm}
\eeq
Using Equation \ref{tau} and Condition \ref{c1} it  is easy to check that the functional $f_{[g]}$ satisfies $f_{[g]} \in C_{\lam}(h)$.  Next it will be shown that 
\bq
f_{[gm_0]}=\frac{\lam(m_0^{-1},h)}{\tau_h(gm_0, m_0^{-1})}f_{[g]}
\eq
for all $m_0 \in M$. This implies that $\{f_{[b_i]}\}_{i=1}^s$ is a basis of $C_{\lam}(h)$. One has:
\beqarn
f_{[gm_0]} & = & \sum_{m \in M}\frac{\lam(m,h)}{\tau_h(gm_0,m)}p_{gm_0m}=
\\ &= & \sum_{n \in M}\frac{\lam(m_0^{-1}n,h)}{\tau_h(gm_0, m_0^{-1}n)}p_{gn}
%\\ &= & \sum_{n \in M}\frac{\lam(g,h)\lam(n,h)\tau_h(g, m)}{\tau_h(gm_0, m_0^{-1}n)}p_{gn} 
\\ &= &\lam(m_0^{-1},h)\sum_{n \in M}\frac{\lam(n,h)}{\tau_h(m_0^{-1},n)\tau_h(g,n)}p_{gn}\frac{\tau_h(g,n)}{\tau_h(m_0^{-1},n)\tau_h(gm_0, m_0^{-1}n)}
\eeqarn

Note that by Equation \ref{tau} one has 
\bq
\frac{\tau_h(g,n)}{\tau_h(m_0^{-1}, n)\tau_h(gm_0, m_0^{-1}n)}=\frac{1}{\tau_h(gm_0, m_0^{-1})}.
\eq
and therefore 
\bq
f_{[gm_0]}=\frac{\lam(m_0^{-1},h)}{\tau_h(gm_0, m_0^{-1})}f_{[g]}
\eq
\epf

\bl
With the above notations one has that
\bq
R_{a, \tau_h}(f_{[g]})=\frac{1}{\tau_h(a, a^{-1}g)}f_{[a^{-1}g]}
\eq
for all $a,g \in G$.
\el
\bpf
Indeed one has that 
\beqarn
R_{a, \tau_h}(f_{[g]}) & = & (f_{[g]}\lh a)\tau_h(a,\;-)\\ & = &( \sum_{m\in M}\frac{\lam(m,h)}{\tau_h(g,m)}p_{a^{-1}gm})\tau_h(a,\;-)
\\ & = &  \sum_{m\in M}\frac{\lam(m,h)}{\tau_h(g,m)}\tau_h(a,\;a^{-1}gm)p_{a^{-1}gm}
\\ & = & \sum_{m\in M}(\frac{\lam(m,h)}{\tau_h(a^{-1}g,m)}p_{a^{-1}gm})(\frac{\tau_h(a^{-1}g,m)\tau_h(a,\;a^{-1}gm)}{\tau_h(g,m)})\eeqarn

Using Equation \ref{tau} observe that 
\bq
\frac{\tau_h(a^{-1}g,m)\tau_h(a,\;a^{-1}gm)}{\tau_h(g,m)}=\frac{1}{\tau_h(a, a^{-1}g)}
\eq
and then the conclusion of the Lemma follows from here.
\epf
Note that for any $u \in k^G$ and $f \in F$ one can write \bq(f.u)\lh m=f.(u \lh (f\rhd m))\eq for all $m \in M$.
%%%%%%%%%%%%%%%%%%%%%%%%%%%%%%%%%%%%%
\bp\label{intord} Let $A\cong k^G\;{^\tau}\dzs kF$ be a cocentral abelian extension of $kF$ by $k^G$. With the following notations it follows that $C(M,\;H,\;\lam)$ is a left coideal subalgebra of $A$.
\ep

\bpf
Suppose that $f\dzs h\in C_{\lam}(h)\dzs h$. Then using Lemma \ref{commf} it follows that
\bq
\Delta(f\dzs h)=\sum_{a\in G}(p_a\dzs h)\ot (R_{a,\tau_h}(f)\dzs h)
\eq
But by Equation \ref{comlr} one has $L_{m,\tau_h}(R_{a, \;\tau_h}(f))=R_{a, \;\tau_h}(L_{m,\tau_h}(f))=\lam(m,h)R_{a, \;h}(f)$ for all $m \in M$ it follows that $R_{a,\tau_h}(f)\in C_{\lam}(h)$ for all $a \in G$. This shows that $\Delta_A(C_{\lam}(h)\dzs h)\subset A \ot (C_{\lam}(h)\dzs h)$ and therefore $C(M,\;H,\;\lam)$ is a left coideal of $A$.

In order to show that $C(M,\;H,\;\lam)$ is an algebra one has to check the following inclusion 
\beq
C_{\lam}(h)(h.C_{\lam}(l))\sigma_{-}(h,l)\subset C_{\lam}(hl),
\eeq
for all $hl \in F$.  For all $f \in C_{\lam}(h)$ and all $g \in C_{\lam}(l)$ one has that
\begin{eqnarray*}
  && R_{m,\tau_{hl}}(f(h.g)\sigma_{-}(h,l)) = [(f(h.g)\sigma_{-}(h,l))\lh m]\tau_{hl}(m,\;-)\\
   &=& (f\lh m)((h.g)\lh m)(\sigma_{-}(h,l)\lh m)\tau_{hl}(m,\;-) \\
   &=&  R_{m,\tau_h}(f)(h.(g\lh (h\rhd m))(\sigma_{-}(h,l)\lh m)\tau_{hl}(m,\;-)\tau^{-1}_{h}(m,\;-)\\
   &=& [\lam(m,h)f](h.[(R_{(h\rhd m)}(g))\tau^{-1}_l((h\rhd m),\;-))\\& &[(\sigma_{-}(h,l)\lh m)\tau_{hl}(m,\;-)\tau_{h}(m,\;-)^{-1}]\\
   &=&  \lam(m,\;h)\lam(h \rhd m,\;l)\sg_m(h,l)f(h.g)\sigma_{-}(h,\l)\\& & \sg^{-1}_m(h,l)\sg^{-1}_{-}(h, l)(\sigma_{-}(h,l)\lh m)\\ & & \tau_l^{-1}((h\rhd m),\;h \rhd \;-)
   \tau_{hl}(m,\;-)\tau_{h}(m,\;-)^{-1}\\
    &=& \lam(m, hl)f(h.g)\sigma_{-}(h,\l)
\end{eqnarray*}
We used that 
\bq
\sg^{-1}_m(h,l)\sg^{-1}_{-}(h, l)(\sigma_{-}(h,l)\lh m)\tau_l^{-1}((h\rhd m),\;h \rhd \;-)=1  
\eq
by Pentagon Equation \ref{sg2-pent} and \bq \lam(m, hl)= \lam(m,\;h)\lam(h \rhd m,\;l)\sg_m(h,l)\eq by Equation \ref{c2}.
Then since $$R_{m,\tau_{hl}}(f(h.g)\sigma_{-}(h,l))=\lam(m, hl)f(h.g)\sigma_{-}(h,\l)$$ it follows by definition of $C_{\lam}(hl)$ that $f(h.g)\sigma_{-}(h,\l)\in C_{\lam}(hl)$.
\epf

%____________________________________________________
\subsection{Left coideal subalgebras of $A$} In this subsection we show that any left coideal subalgebra of $A$ is of the type $C(M,\;H,\; \lam)$ as above.
This result is inspired by \cite{Grln} where the left coideal subalgebras where described under the additional assumption of a $"^*"$-structure on $A$. We will give in this subsection a proof that does not use this assumption but  it uses the structure of coideal subalgebras from Theorem \ref{dec} instead.

With the notations from the previous Subsection define the following linear functionals $\delta_{b_iM}:=\sum_{m\in M}p_{b_im} \in k^G$ for all $i =1, s$. 
%We will prove, Theorem \ref{main}, namely that any left coideal subalgebra of $A$ is of the type $C(M,\;H,\;\lam)$ defined above.

{\bf Proof of Theorem \ref{main}:}
\md
Let $B$ be an arbitrary left coideal subalgebra of $A$. Write the elements of $B$ as $b=\sum_{h\in F}f_h\#h$ where $f_h \in k^G$. 
Note that by Lemma \ref{commf} 
$$
\Delta_A(b)=\Delta_A(\sum_{h\in F}f_h\#h)=\sum_{h\in F,\\ \;g\in G}(p_g\#h)\otimes (R_{g, \tau_h}(f_h)\#h).
$$
Since $B$ is a left coideal, it follows that for each fixed pair of elements $(g, \;h) \in G\times F$ one has that $R_{g, \tau_h}(f_h)\# h \in B$.

Therefore one can write $B = \oplus_{h\in F}(B(h)\#h)$ with $B(h)$ a subspace of $k^G$ which is mapped by $R_{g, \tau_h}$ to itself for all $g\in G$. Since $B$ is also an algebra, we have
$$B(h_1)(h_1.B(h_2))\sg(h_1,h_2) \subseteq B(h_1h_2),$$ for all $h_1, h_2 \in F$.

In particular this implies that $B(1)$ is a subalgebra of $k^G$ which affords a left representation of $G$ via the operators $R_{g, \tau_1}(f)=f \lh g$. It follows from Lemma \ref{existence} that there is a subgroup $M \leq G$ such that $B(1)$ is the subspace $k^{(G/M)_l}$ of $M$-left invariant functions on $G$. Let $G = \sqcup_{i} b_iM$, $1 \leq  i \leq r$ be the left coset decomposition of $G$ with respect to the subgroup $M$. Here $r = |G|/|M|$ is the index of $M$ in $G$. Then the linear functionals $\delta_{b_iM}$  form a linear basis on the algebra $B(1)$.

Let $H := \{h \in F|B(h) \neq  0\}$. It will be shown next that $H$ is a subgroup of $F$. Let $y_B$ be the left invariant element of $B$ defined in Subsection \ref{integral}. Suppose further that $$y_B=\sum_{h \in H'}u_h\dzs h\in B$$ for some nonzero elements $u_h\in B(h)$ and some non-empty subset $H'\subset H$.
%Since $y_B^2=y_B$ it follows.

By Theorem \ref{dec} since $B=y_B\leftharpoonup A^*$ it follows that $H'=H$. Moreover the same Theorem implies that $B(h)=<\{L_{a, \tau_h}(u_h)\}_{a \in G}>$,  the linear span  of the functionals $\{L_{a, \tau_h}(u_h)\}_{a \in G}$.  Indeed, $A^*$ ca be identified as algebras to $k^F\#_{\tau^*} kG$ via 
\beq
<p_x \#_{\tau^*}a, p_b\dzs y>=\delta_{a,b}\delta_{x,y}
\eeq
\ncm{\dzt}{\#_{\tau^*}}
for all $x, y\in F$ and $a,b \in G$.
Then using Equation \ref{commf} one has
\beqarn
(u_h \dzs h) \lh (p_x \dzt a)& = & \sum_{g \in G}<p_x \dzt a, p_g\dzs h>R_{g ,\tau_h}(u_h)\dzs h\\ &= & \delta_{x,h}R_{a ,\tau_h}(u_h).\eeqarn Thus $y_B \lh p_h\dzt a= R_{a ,\tau_h}(u_h)\dzs h$. 
\md
Next it will be shown that $H$ is a subgroup of $F$. Since $(\delta_{M}\#1)y_B=y_B$ it follows that $\delta_Mu_h=u_h$ for all $h \in H$. Since $u_h\neq 0$ it follows that $\mtr{supp}(u_h)\subseteq M$. Thus one has $u_h\in B(h)\delta_{M}$ for all $h\in H$. In particular $u_1\in B(1)\di_M=k\delta_M$.

Without loss of generality one may suppose further that $u_1=\delta_{M}$. %Then clearly $y_B^2=\al y_B$ for some nonzero scalar $\al \in k^*$.
Since $B(h)\neq 0$ it follows that there is $u \in B(h)$ and $x \in G$ such that $u(x)\neq 0$. It follows that $R_{x ,\tau_h}(u)(1)=u(x)\neq 0$. On the other hand since $R_{x, \tau_h(f)}\in B(h)$ one can conclude that there is $f:=R_{x ,\tau_h}(u) \in B(h)$ with $f(1)\neq 0$. Therefore for such element $f \in k^G$ one has that $(f\dzs h)y_B=f(1)y_B$ is a nonzero element of $B$. On the other hand since 
\bq\lb{f1}
(f\dzs h)y_B=\sum_{l \in H}f(h.u_l)\sg(h,l)\dzs hl
\eq
we deduce that $hH\subseteq H$. Thus $H$ is a subgroup of $F$. Moreover from Equation \ref{f1} we deduce that
\beq
f(h.u_l)\sg(h,l)=f(1)u_{hl}
\eeq
for all $l \in H$. For $l=h^{-1}$ this identity becomes
\beq\lb{idf}
f(h.u_{h^{-1}})\sg(h,h^{-1})=f(1)\delta_M
\eeq

Evaluating both sides of the last identity at $g=1$ it follows that $u_h(1)=1$ for all $h \in H$. Moreover since $(\delta_M\dzs 1)y_B=y_B$ it follows that $\delta_Mu_h=u_h$ for any $h\in H$. This shows that $\mtr{supp}(u_h)\subseteq M$. On the other hand Equation \ref{idf} shows that $M \subseteq \mtr{supp}(f)$ for any $f \in B(h)$ with $f(1)\neq 0$. In particular, since $u_h(1)=1$ it follows that $\mtr{supp}(u_h)=M$. But evaluating Equation \ref{idf} at any $m \in M$ it follows that:
\beq
f(m)u_{h^{-1}}(h \rhd m)\sg_m(h, h^{-1})=f(1)
\eeq
 Since $\mtr{supp}(u_h)=M$ this shows that $M$ is stable under the action of $H$ and $\dim B(h)\delta_M$=1. Thus  $B(h)\delta_M=<u_h>$.

Since $$(u_h\dzs h)y_B=y_B$$ this implies that
\begin{equation}\label{multica}
    u_h(h.u_{h'})\sg(h,h')= u_{hh'}
\end{equation}
for all $h,h'\in H$. Evaluating both sides of the last equality at $m\in M$ one has
\begin{equation}\label{multicaev}
    u_h(m)u_{h'}(h\rhd m)\sg_m(h.h')=u_{hh'}(m)
\end{equation}
%%%%
Since for $m \in M$ the operator $L_{m, \tau_h}$ maps $B(h)\delta_{M}$ to itself it follows that there is a function $\lam : M \times H \ra k^*$ such that
\beq\lb{ondil}
L_{m, \tau_h}(u_h) = \lam(m, h)u_h.
\eeq
Evaluating  both sides of the above identity at $1$ one gets that $u_h(m)=\lam(m,\;h)$ for all $m\in M$. On the other hand evaluating both sides of the same equality at $m' \in M$ it follows that $u_h(mm')\tau_h(m,m')=\lam(m,\;h)u_h(m')$ and therefore one obtains Equation \ref{c1}, i.e, 
\beq
\lam(mm',\; h)\tau_h(m,m')=\lam(m,\;h)\lam(m,\;h')
\eeq for all $m,m'\in M$ and $h \in H$.
For the other Equality \ref{c2} note that 
\begin{eqnarray*}
  \lam(m,\;hh') &=& u_{hh'}(m)\\
   &=& u_h(m)u_{h'}( h\rhd m)\sg_m(h,h') \\
   &=& \lam(m,\;h)\lam( h\rhd m, \; h')\sg_m(h,h')  \\
\end{eqnarray*}
It remains to show that $B(h)$ coincides to $C_{\lam}(h)$, the subspace of all functions $f\in k^G$ which verify $$L_{m, \tau_h}(f) = \lam(m, h)f,$$ for all $m \in M$. Note that since $B(h)B(1) \subseteq B(h)$ it follows that $B(h) = \oplus_{i=1}^sB(h)\delta_{b_iM}$. Since $L_{b_i, \tau_h}(B(h)\di_M)=B(h)\di_{b_iM}$ it follows that all the spaces $B(h)\delta_{b_iM}$ are also one dimensional. Thus $\mtr{dim}_k(B(h))=\frac{|G|}{|M|}$. On the other hand Lemma \ref{desciptch} shows that $\mtr{dim}_k(C_{\lam}(h))=\frac{|G|}{|M|}$.

Thus in order to conclude that $B(h)=C_{\lam}(h)$ it is enough to show that $L_{m,\tau_h}(f) = \lam(m, h)f$ for all $f  \in B(h)$. If $f \in B(h)\delta_M$ then this relation is satisfied from the definition of $\lam$ and Relation \ref{ondil}. On the other hand since $B(h)$ is spanned as vector spaces by $R_{g, \tau_h}(u_h) $ Equation \ref{comlr} implies that for all $g \in G$ one has
\bq
L_{m, \tau_h}(R_{g , \tau_h}(u_h))= R_{g, \tau_h}(L_{m, \tau_h}(u_h))=\lam(m,\;h)R_{g , \tau_h}(u_h).
\eq $\square$. \subsection{On the Wall's conjecture}\lb{wall}
%In the paper \cite{Grln} the authors proposed the following version of  Wall's onjecture:
%{\bf Conjecture 1.:} Let $\mtc{F}$ be a finite dimensional semisimple fusion algebra with $n$ simple objects. Then the number
%of  minimal (maximal) fusion subalgebras which are generated by a subset of the simple objects of $\mtc{F}$ is less than $n$.
%\blue{ maximal coideal subalegrbas ok; minimal coideal subalgebras = maximal normal Hopf subalgebras of $A^*$}
%If we take F to be the group algebra of G, then Conjecture 1.2 is equivalent to Wall�s conjecture.
%Recall that the Grothendieck ring of any finite tensor category is a fusion algebra.  Moreover, following \cite{gmj} every minimal normal fusion subcategory of $\mtr{Rep}(A)$ is given by $\mtr{LKer}_A(M)$ for an irreducible $A$-module $M$.
In the spirit of \cite{Grln} we have the following Theorem which is a Hopf algebra analogue of Wall's conjecture. 
\bt
Let $A\cong k^G\;^{\tau}\#_{\sg} kF$ be a cocentral extension of Hopf algebras with $G$ and $F$ solvable groups. Then the number of maximal (resp. minimal) left coideal subalgebras of $A$ is less or equal than the dimension of $A$.
\et
With the above characterization of coideal subalgebras the proof of this Theorem is the same as that of Theorem 3.8 from \cite{Grln}.

\br
In the same paper \cite{Grln}, the authors proposed Conjecture 1.2, as a new version of Wall's conjecture for semisimple fusion algebras. It was announced that this conjecture is solved  for commutative fusion rings. We remark that there are cocentral Kac algebras with noncommutative Grothendieck rings, for example  the smash product Hopf algebra $k^{Q_8} \# kC_2$, dual to the smash coproduct Hopf algebra from Section 8 of \cite{k16}.
\er

\section{Hopf subalgebras of cocentral abelian extensions of Hopf algebras}
In this section we describe all Hopf subalgebras of an abelian cocentral extension.
\bt\label{hsa} Let $A=k^G\;^{\tau}\#_{\sg} kF$ be a cocentral abelian extension of $kF$ by $k^G$. Then all Hopf subalgebras of $A$ are of the form $C(M,\;H,\;\lam)$ where $M$ is a normal subgroup of $G$ and $\lam$ satisfies the additional invariance condition:
\bq\lb{hl}
\lam(a^{-1}ma,h)=\lam(m,h)\tau_h^{-1}(a, a^{-1}ma)\tau_h(m,a)
\eq
for all $a \in G$ and $m \in M$.
\et
\bpf
A Hopf subalgebra $B$ of $A$ is in particular a left coideal subalgebra and by the previous Theorem is of the type $C(M,\;H,\;\lam)$ defined above. Suppose that $B:=C(M,\;H,\;\lam)$ is a Hopf subalgebra of $A$ and let $f \#_{\sg} h \in B$. Then using Lemma \ref{commf} one has
$$
\Delta_A(f \#_{\sg} h)= \sum_{a \in G}(L_{a,\tau_h}(f)\dzs h)\ot (p_a\#_{\sg} h)%=\sum_{a \in G}  ((\sum_{u  \in G}f_a\tau_h(au^{-1}, u)p_{au^{-1}})\#_{\sg} h)\otimes(p_u \#_{\sg} h).
$$ Then  $\Delta_A(f \#_{\sg} h)\in B\otimes A$ if and only if 
\bq
L_{a, \tau_h}(f)=\sum_{x  \in G}f_x\tau_h(xa^{-1}, a)p_{xa^{-1}}\in C_{\lam}(h)
\eq
for any $a\in G$. This implies that for any $a \in G$ the space $Ma^{-1}$ is also a right coset of $M$ in $G$, i.e $M$ is a normal subgroup of $G$. On the other hand using repeatedly Equation \ref{leff} it follows that:
\beqarn
L_{m, \tau_h}(L_{a,\tau_h}(f))&=&L_{ma,h}(f)\tau_h(m,a)\\ & =& L_a(L_{a^{-1}ma, \tau_h}(f))\tau_h^{-1}(a, a^{-1}ma)\tau_h(m,a)\\ &=& \lam(a^{-1}ma,h)\tau_h^{-1}(a, a^{-1}ma)\tau_h(m,a)L_{a,\tau_h}(f)% \lam(m,h)L_{a, \tau_h}(f)
\eeqarn
This shows that $L_{a, \tau_h}(f)\in C_{\lam}(h)$ if and only if \bq\lam(a^{-1}ma,h)=\lam(m,h)\tau_h^{-1}(a, a^{-1}ma)\tau_h(m,a).\eq
\epf
\br

Alternatively, for the converse of the above Theorem it can be shown that
\begin{equation*}
   \Delta_{k^G}(f_{[a]})=\sum_{y \in G}\frac{\lam(y^{-1},h)}{\tau_h(a,y^{-1})}f_{[yay^{-1}]}\ot p_{y}
\end{equation*}
and thus $\Delta_{k^G}(C_{\lam}(h))\subset C_{\lam}(h)\ot  k^G$.
\er
%%%%%%%%%%%%%%%%%%%%%%%%%%%%%%%%%%%%%%%%%%%%%%%%%%%%%
\bibliographystyle{amsplain}
\bibliography{macbob}
\end{document}

%% file: macros.tex
\ncm{\cE}{\mtc{E}}
\ncm{\End}{\mtr{End}}
\ncm{\cop}{\mtr{cop}}
\ncm{\op}{\mtr{op}}
\ncm{\chr}{character }
\ncm{\bw}{\bwt}
\ncm{\cY}{\mtc{Y}}
\ncm{\bx}{\boxtimes}

\ncm{\eeq}{\end{equation}}

\ncm{\bea}{\begin{eqnarray}}
\ncm{\eea}{\end{eqnarray}}
\ncm{\beanon}{\begin{eqnarray*}}
\ncm{\eeanon}{\end{eqnarray*}}\ncm{\ek}{\eps|_K}\ncm{\diez}{\#}
\ncm{\bwt}{\bowtie}
\ncm{\cC}{\mtc{C}}